\documentclass[12pt]{iopart}
\usepackage{amssymb}
\usepackage{amsthm,harvard}
\newtheorem{thm}{Theorem}
\newtheorem{conjecture}{Conjecture}
\newtheorem{cor}{Corollary}
\newtheorem{lem}{Lemma}
\makeatother

\begin{document}

\title[On some results of Cufaro Petroni about Student t-processes]{On some results of Cufaro Petroni about Student t-processes}

\author{C. Berg$^1$ and C. Vignat$^2$}

\address{$^1$ Department of Mathematics, University of Copenhagen, Universitetsparken 5, DK-2100, Copenhagen, Denmark\\
$^2$ Institut Gaspard Monge, Universit\'{e} de Marne la 
Vall\'{e}e, France}
\eads{\mailto{berg@math.ku.dk}, \mailto{vignat@univ-mlv.fr}}

\begin{abstract}
This paper deals with Student t-processes as studied in \cite{cufaro}. We prove and extend some conjectures expressed by Cufaro Petroni about the asymptotical behavior of a Student t-process and the expansion of its density. First, the explicit asymptotic behavior of any real positive convolution  power of a Student t-density with any real positive degrees of freedom is given in the multivariate case; then  the integer convolution power of a Student t-distribution with odd degrees of freedom is shown to be a convex combination of Student t-densities with odd degrees of freedom. At last, we show that this result does not extend to the case of non-integer convolution powers.
\end{abstract}
\maketitle

\section{Introduction}

In a recent contribution \cite{cufaro}, N. Cufaro Petroni derived several
results about the behavior of some non stable L\'{e}vy processes with Student t-marginals. 
More precisely, he considered the random walk
\[
Z_{N}=\sum_{i=1}^{N}X_{i}
\]
where $N\in\mathbb{N}$ and each independent step $X_{i}$ follows
a Student t-distribution with $f=2n+1$ degrees of freedom, $n\in\mathbb{N}$.
We recall that the Student t-density with $f=2\nu$ degrees of freedom ($\nu>0$) is
\[
f_{\nu}\left(x\right)=A_{\nu}\left(1+x^{2}\right)^{-\left(\nu+\frac{1}{2}\right)}; \,\,\, A_{\nu}=\frac{\Gamma\left(\nu+\frac{1}{2}\right)}{\Gamma\left(\frac{1}{2}\right)\Gamma\left(\nu\right)}.
\]

The family of Student t-densities includes the Cauchy density for
$f=1$ and the scaled density
$\frac{1}{\sqrt{\nu}}f_{\nu}(\frac{x}{\sqrt{\nu}})$ converges to the
Gaussian density  as  $f\rightarrow + \infty.$ All Student
t-distributions are heavy tailed. Grosswald \cite{grosswald}
proved that they are infinitely divisible. They also have the stronger
property of being self-decomposable, cf. \cite{steutel}.

Stochastic processes with Student t-marginals and various types of dependence structures have been proposed in \cite{Heyde}, most of them with dependent increments. On the other side, Cufaro Petroni's paper deals with L\'{e}vy Student t-processes, which exist by the infinite divisiblity of the Student t-distribution. In both cases, these processes have  applications in finance \cite{Schoutens} and in physics \cite{Vivoli}.
\section{Three conjectures by Cufaro Petroni}

Cufaro Petroni obtained precise results about the process $Z_{N}$ only in the case
of $f=3$ degrees of freedom; however, he expressed three
conjectures about the extension of these results to more general cases; 
the first conjecture is

\begin{conjecture}
\label{con:conj1}
For all $N\in\mathbb{N}$ and for all $f=2n+1,\,\, n\in\mathbb{N}$,
the distribution of $N^{-1}Z_{N}$ is a convex combination of Student
t-distributions with odd degrees of freedom.
\end{conjecture}
The two remaining conjectures concern the distribution of $Z_{N}$
for non-integer values of $N$, which makes sense because the Student t-distribution is infinitely divisible: the $c-$fold convolution of the distribution
$f_{\nu}$ is defined, for any real positive $c,$ as the inverse
Fourier transform\[
f_{\nu}^{*c}\left(x\right)=\frac{1}{2\pi}\int_{-\infty}^{+\infty}e^{iux}\left[\varphi_{\nu}\left(u\right)\right]^{c}du,\]
where $\varphi_{\nu}\left(u\right)$ is the characteristic function
of the Student t-distribution\[
\varphi_{\nu}\left(u\right)=k_{\nu}\left(\vert u\vert\right)\]
with 
\begin{equation}
\label{knu}
k_{\nu}\left(u\right)=\frac{2^{1-\nu}}{\Gamma\left(\nu\right)}u^{\nu}K_{\nu}\left(u\right),  u>0.
\end{equation}
Here $K_{\nu}$ is the modified Bessel function of the second kind
also called the Macdonald function. 
The expression (\ref{knu}) reduces to elementary functions exactly
when $\nu=n+1/2,n=0,1,\ldots$ because
\begin{equation}
\label{knuodd}
k_{n+\frac12}(u)=e^{-u}q_n(u),  u>0,
\end{equation}
where $q_n$ is a polynomial of degree $n$ with positive coefficients,
called the $n$'th Bessel polynomial. It is given as
\begin{equation}\label{eq:qn}
q_{n}\left(u\right)=\sum_{k=0}^{n}\alpha_{k}^{(n)}u^{k},
\end{equation}
where 
\begin{equation}\label{eq:alpha}
\alpha_{k}^{(n)}=\frac{(-n)_{k}2^k}{(-2n)_{k}k!}.
 \end{equation}
The first examples of these polynomials
are
\[
q_{0}\left(u\right)=1, q_{1}\left(u\right)=1+u, q_{2}\left(u\right)=1+u+\frac{u^{2}}{3},
\]
cf.  \cite{Berg}.

Cufaro Petroni's second conjecture concerns the asymptotic behavior of the density of the $c-$fold convolution of $f_{\nu}$:

\begin{conjecture}
\label{con:conj2}For every $c>0$ and $\nu>0$, the asymptotic behavior
of the $c-$fold convolution $f_{\nu}^{*c}$ is given by
\[
f_{\nu}^{*c}\left(x\right)\sim\frac{cA_{\nu}}{x^{2\nu+1}}, \,\, x  \rightarrow+\infty.
\]
\end{conjecture}
Cufaro Petroni's last conjecture is an extension of Conjecture 1 to the $c-$fold convolution $f_{\nu}^{*c}$ as follows:
 
\begin{conjecture}
\label{con:conj3}
Conjecture \ref{con:conj1} extends to non-integer sampling
times $c$ under the following form: for all $\nu_{0}>0$ and all $c>0$, 
\[
f_{\nu_{0}}^{*c}(x) = \int_{\nu_{0}}^{+\infty} q_{\nu_{0},c}(\nu) \frac{1}{c}f_{\nu}(\frac{x}{c}) dx
\]
for some distribution $q_{\nu_{0},c}(\nu)$. 
\end{conjecture}

In this paper, we show that conjecture (\ref{con:conj1}) holds true
and give an extended version of it; likewise, we prove the conjecture
(\ref{con:conj2}). We were unable to prove or disprove conjecture (\ref{con:conj3}), but we disprove a discrete version of it in the case where $\nu=n+\frac{1}{2}$ with $n\in\mathbb{N}.$

Moreover, we consider in the rest of this paper the multivariate context: all Student t-variables are supposed rotation invariant $d-$dimensional vectors. The multivariate Student t-density is given, for $\mathbf{x}=\left( x_{1},\dots,x_{d} \right) \in \mathbb{R}^{d}$ by

\[
f_{\nu}(\mathbf{x}) = A_{d,\nu} \left( 1+\vert \mathbf{x}\vert ^2 \right)^{-\nu-d/2}, \,\,\, A_{d,\nu}=\frac{\Gamma(\nu+\frac{d}{2})}{\Gamma(\nu) \Gamma(\frac{1}{2})^d},
\]

where 

\[
\langle \mathbf{x},\mathbf{y}\rangle = \sum_{i=1}^{d} x_{i}y_{i},\,\,\, \vert \mathbf{x} \vert= \langle \mathbf{x},\mathbf{x}\rangle ^{\frac{1}{2}},\,\,\, \mathbf{x},\mathbf{y} \in \mathbb{R}^{d}.
\]

\section{First conjecture: the odd degrees of freedom case}

Cufaro Petroni's first conjecture \cite[Prop. 5.2]{cufaro} is that if $X_{i}$
is a set of independent Student t-distributed random variables with
$f=2n+1, n\in\mathbb{N}$ degrees of freedom, then the density of the distribution
of the normalized $N-$th step of the random walk\[
N^{-1}Z_{N}=\frac{1}{N}\sum_{i=1}^{N}X_{i}\]
writes as \begin{equation}
v\left(x\right)=\sum_{k=n}^{nN}\beta_{k}^{\left(n,N\right)}f_{k+\frac{1}{2}}\left(x\right)\label{eq:fz}\end{equation}
with $\beta_{k}^{\left(n,N\right)}\ge0,n\le k\le nN.$ We extend
and prove this conjecture as follows:

\begin{thm}
\label{thm:thm1}
If $N\in\mathbb{N}$ and 
\[
\mathbf{Y}_{N}=\sum_{i=1}^{N}a_{i}\mathbf{X}_{i},
\]
where $a_{i}$ are positive numbers with sum $1$ and $\mathbf{X}_{i}$ are
independent $d-$variate Student t-distributed, each with $f_{i}=2n_{i}+1 \left(n_{i}\in\mathbb{N}\right)$
degrees of freedom, then the density of $\mathbf{Y}_{N}$ is \[
\sum_{j=min\left(n_{1},\dots, n_{N}\right)}^{n_{1}+\dots+n_{N}}\beta_{j}f_{j+\frac{1}{2}}\left(\mathbf{x}\right)\]
where the coefficients $\beta_{j}$ are nonnegative with sum $1$ and do not depend on the dimension $d.$
\end{thm}
\begin{proof}
The characteristic function of the $d-$variate Student t-distribution is
\[
\varphi\left( \mathbf{u}\right) = k_{\nu}(\vert \mathbf{u} \vert)
\]
where the function $k_{\nu}$ is given by (\ref{knu}). Since  for $\nu_{i} = n_{i}+\frac{1}{2}$, this function reads  $k_{\nu_{i}}(\vert \mathbf{u} \vert) = e^{-\vert \mathbf{u} \vert}q_{n_{i}}(\vert \mathbf{u} \vert)$ where $q_{n_{i}}$ is the Bessel polynomial of degree $n_{i}$, the result follows from  \cite[Th. 2.6]{Berg}:
\begin{equation}\label{eq:final}
q_{n_1}(a_1u)q_{n_2}(a_2u)\cdots
q_{n_N}(a_Nu)=\sum_{j=l}^L\beta_jq_j(u),\quad u\in\mathbb R
\end{equation}
with nonnegative coefficients $\beta_j$ with sum 1 and 
$l=\min(n_1,\ldots,n_N),L=n_1+\cdots +n_N$.
\end{proof}

As a particular case, choosing $a_{i}=\frac{1}{N}, 1\le i\le N$
for $N\in\mathbb{N}$, we deduce that the coefficients $\beta_{k}^{\left(n,N\right)}$
in (\ref{eq:fz}) are positive, and thus the density $g$ is a convex
combination of Student t-distributions with odd degrees of freedom.

We are not able to provide an expression for the coefficients
$\beta_{k}^{\left(n,N\right)}$ which can be used directly to see the
non-negativity. Using Carlitz' formula, see \cite{Berg},
\begin{equation}\label{eq:inverse}
u^n=\sum_{j=0}^n \delta_j^{(n)}q_j(u),\quad n=0,1,\ldots
\end{equation}
with 
\begin{equation}\label{eq:carlitz}
 \delta_{j}^{\left(n\right)}=\left\{ \begin{array}{cc}
\frac{\left(n+1\right)!}{2^{n}}\frac{\left(-1\right)^{n-j}\left(2j\right)!}{\left(n-j\right)!j!\left(2j+1-n\right)!}
& 
\mbox{ for }\frac{n-1}{2}\le j\le n\\
0 & \mbox{ for } 0\le j<\frac{n-1}{2}\end{array}\right.
\end{equation}
it is possible to write
$$
\prod_{j=1}^N
q_{n_j}(a_ju)=\sum_{k_1=0}^{n_1}\cdots\sum_{k_N=0}^{n_N}\left(\prod_{j=1}^N \alpha_{k_j}^{(n_j)}
a_j^{k_j}\right)\sum_{i=0}^{k_1+\cdots
  +k_N}\delta_i^{(k_1+\cdots+k_N)}q_i(u),
$$
which gives an expression for $\beta_{j}$ in (\ref{eq:final}), but because
of the varying sign of $\delta_j^{(n)}$, it is not possible to see
directly that $\beta_j\ge 0$.

If $a_j=1/N$ and $n_1=\ldots=n_N=1$, i.e. the case of $f=3$ degrees of
freedom where $q_1(u)=1+u$, this formula simplifies to the expression
given in \cite[Prop. 5.2]{cufaro}. It is claimed that the expression
is positive, but no convincing argument is given.

\section{Second conjecture: the asymptotic behavior of the Student process}

A second property studied by Cufaro Petroni is the asymptotic behavior of
the distribution of the random walk $Z_{N}$; in the case of $f=3$
degrees of freedom $\left(\nu=\frac{3}{2}\right),$ he obtains
the following result \cite[Prop. 5.1]{cufaro}: for all $c>0,$\[
f_{\frac{3}{2}}^{*c}\left(x\right)\sim\frac{2c}{\pi x^{4}}, \,\, x \rightarrow+\infty.\]
We provide now an extension of this result to any value $f=2\nu \,\, (\nu>0)$
of degrees of freedom. Cufaro Petroni's argument is via Fourier analysis. This argument becomes very technical if one tries to generalize it to arbitrary degrees of freedom. Our proof is based on results about subexponential distributions.

\begin{thm}
\label{thm:cstepStudent}
For any $c>0$ and any $\nu>0,$ the density of the $c-$fold
convolution of the $d-$variate Student t-distribution behaves asymptotically as
\[
f_{\nu}^{*c}\left(\mathbf{x} \right)\sim\frac{cA_{d,\nu}}{\vert \mathbf{x} \vert^{2\nu+d}}, \,\, \vert \mathbf{x} \vert \rightarrow +\infty.
\]
\end{thm}
\begin{proof}
The proof is based on a series of lemmas given in the last section.
The $d-$variate Student t-distribution is subordinated to the $d-$variate Gaussian semigroup
\[
g_{t}(\mathbf{x} )=(4\pi t)^{-d/2} \exp{\left( -\frac{\vert \mathbf{x} \vert^{2}}{4t} \right)},
\,\,
t>0, \,\, \mathbf{x}\in \mathbb{R}^{d}
\]
by the inverse Gamma density, i.e.
\[
f_{\nu}\left(\mathbf{x}\right)=\int_{0}^{+\infty}g_{t}\left(\mathbf{x}\right)dH_{\nu}\left(t\right),\]
where $H_{\nu}\left(t\right)$ is the inverse Gamma distribution
with density
\[
h_{\nu}\left(t\right)=C_{\nu}\exp\left(-\frac{1}{4t}\right)t^{-\nu-1},\,\, t>0,\,\, C_{\nu}=\frac{1}{2^{2\nu}\Gamma\left(\nu\right)}.
\] 
From this representation
we deduce in Lemma \ref{lem:cconvolutionpower} the same representation
for the c-fold convolution power of the Student t-density, namely
\[
f_{\nu}^{*c}\left(\mathbf{x}\right)=\int_{0}^{+\infty}g_{t}\left(\mathbf{x}\right)dH_{\nu}^{*c}\left(t\right).
\]

We note that this property is very general, in the sense that it holds
for any infinitely divisible probability distribution $dH\left(t\right)$
on $\left[0, \infty \right[.$ 

The next step of the proof is the derivation of the asymptotic behavior of
the c-fold convolution power of the inverse Gamma density
$h_{\nu}\left(t\right)$: by Lemma \ref{lem:inverseGammaasymptotic},
this reads\[
h_{\nu}^{*c}\left(t\right)\sim cC_{\nu}t^{-\nu-1},\,\,\, t\rightarrow+\infty.\]

Finally, we show in Lemma \ref{lem:subordination} that this asymptotic
behavior implies, by subordination to the Gaussian semigroup, the
desired asymptotic behavior of the c-fold Student t-convolution.
\end{proof}
As a consequence of this theorem, we deduce the following

\begin{cor}
In the case where the number of degrees of freedom 
$2\nu =2n+1$ is an odd integer and with integer $N,$ the coefficient $\beta_{n}^{\left(n,N\right)}$
in (\ref{eq:fz}) reads \[
\beta_{n}^{\left(n,N\right)}=\frac{1}{N^{2n}}.\]
\end{cor}
\begin{proof}
Since the coefficients $\beta_{k}^{(n,N)}$ do not depend on the dimension $d$, we consider the case $d=1.$
The function $v$ in (\ref{eq:fz}) is the density of the normalized random
walk $\frac{1}{N}\sum_{i=1}^{N}X_{i}$ and thus writes
\[
v\left(x\right)=Nf_{n+\frac{1}{2}}^{*N}\left(Nx\right).
\]
By Theorem \ref{thm:cstepStudent}, 
\[
v\left(x\right)\sim N^{2}\frac{A_{n+\frac{1}{2}}}{\left(Nx\right)^{2n+2}}=\frac{A_{n+\frac{1}{2}}}{N^{2n}}x^{-2n-2}, \,\, x \to + \infty.\]
Since each Student t-distribution $f_{k+\frac{1}{2}}$ in (\ref{eq:fz})
has asymptotic behavior
\[
f_{k+\frac{1}{2}}\left(x\right)\sim A_{k+\frac{1}{2}}x^{-2k-2},\]
we deduce that
\[
v\left(x\right)\sim\beta_{n}^{\left(n,N\right)}A_{n+\frac{1}{2}}x^{-2n-2}.
\]
Identification of the two equivalents yields the result.
\end{proof}

\section{Third conjecture: non integer sampling time and odd degrees of freedom}

In this section, we prove by contradiction the following result

\begin{thm}
\label{thm:thm3}
For all $c>0,\,\, c\notin\mathbb{N}$ and $\nu = n+\frac{1}{2},\,\,n\in\mathbb{N},$ the $1-$variate density
$f_{\nu}^{*c}$ can not be expanded as
\begin{equation}
\label{fnu}
f_{n+\frac{1}{2}}^{*c}\left(x\right)=\sum_{j=0}^{+\infty}\beta_{j}\frac{1}{c}f_{j+\frac{1}{2}}\left(\frac{x}{c}\right),
\end{equation}
with parameters $\beta_{j} \ge 0.$ 
\end{thm}
\begin{proof}
We remark that integrating equality (\ref{fnu}) over $\mathbb{R}$ yields $\sum_{j=0}^{+\infty} \beta_{j}=1$ so that the sequence $(\beta_k)$ is summable.
The  Fourier transform  of (\ref{fnu}) reads
\[
k_{n+\frac{1}{2}}^{c}(u) = \sum_{j=0}^{+\infty} \beta_{j} \exp(-cu) q_{j}(cu),\,\,u>0
\]
where $q_{j}$ is the Bessel polynomial of degree $j.$ Thus, by Lemma \ref{lem:seriesS}, the sum $\sum_{j=0}^{+\infty} \beta_{j} q_{j}(cu)$ is an entire function, so that the function $k_{n+\frac{1}{2}}^c(u)$ extends to an entire function. But 
\[
k_{n+\frac{1}{2}}^{c}(u) = \exp(-cu) \left[ q_{n} (u) \right]^{c}
,\] and since $c$ is not an integer, the function $q_{n}^{c}$ is not holomorphic at any of the complex roots of $q_{n}$, what concludes the proof.
\end{proof}

\section{Lemmas for the proof of theorems \ref{thm:cstepStudent} and \ref{thm:thm3}}

\begin{lem}
\label{lem:cconvolutionpower}
The $c-$fold convolution of the density
$f_{\nu}$ reads, for all $c>0$ and $\nu>0,$\[
f_{\nu}^{*c}\left( \mathbf{x} \right)=\int_{0}^{+\infty}g_{t}\left(\mathbf{x}\right)dH_{\nu}^{*c}\left(t\right).\]

\end{lem}
\begin{proof}
Let us consider the $d-$dimensional Fourier transform $\mathcal{F}$ of the right hand side
\begin{eqnarray*}
\mathcal{F}\left[\int_{0}^{+\infty}g_{t}\left(\mathbf{x}\right)dH_{\nu}^{*c}\left(t\right)\right]\left(\mathbf{y}\right) & = & \left[\int_{0}^{+\infty}\mathcal{F}\left[g_{t}\left(\mathbf{x}\right)\right]dH_{\nu}^{*c}\left(t\right)\right]\left(\mathbf{y}\right)\\
 & = & \int_{0}^{+\infty}\exp\left(-t\vert \mathbf{y}\vert^{2}\right)dH_{\nu}^{*c}\left(t\right).
\end{eqnarray*}
This last integral is nothing but the Laplace transform $\mathcal{L}$ of $H_{\nu}^{*c}$
evaluated at $\vert \mathbf{y} \vert^{2}$, and thus coincides with \[
\mathcal{L}\left(H_{\nu}\right)^{c}\left(\vert \mathbf{y} \vert^{2}\right)=\left(\int_{0}^{+\infty}\exp\left(-t\vert \mathbf{y} \vert^{2}\right)dH_{\nu}\left(t\right)\right)^{c}=\left(\mathcal{F}\left(f_{\nu}\left(\mathbf{x}\right)\right)\right)^{c}\left(\mathbf{y}\right).\]
The result follows by considering the inverse Fourier transform.
\end{proof}

\begin{lem}
\label{lem:inverseGammaasymptotic}
For all $c>0$ and $\nu>0,$ the $c-$fold convolution of the inverse Gamma density has asymptotical behavior 
\[
h_{\nu}^{*c}\left(t\right)\sim cC_{\nu}t^{-\nu-1},\,\,\, t\rightarrow+\infty.
\]
\end{lem}

\begin{proof}
Since\[
h_{\nu}\left(t\right)=C_{\nu}\exp\left(-\frac{1}{4t}\right)t^{-\nu-1}\sim C_{\nu}t^{-\nu-1},\,\,t\rightarrow +\infty,\]
the tail function $\bar{H}_{\nu}(t)=1-H_{\nu}(t)$ of the inverse Gamma distribution has the asymptotic
behavior\[
\bar{H}_{\nu}\left(t\right)\sim\frac{C_{\nu}}{\nu}t^{-\nu},\,\,t\rightarrow
+\infty .\]
This tail function is thus regularly varying and, by \cite[p.278]{Feller},
\[
\bar{H}_{\nu}^{*2}\left(t\right)\sim\frac{2C_{\nu}}{\nu}t^{-\nu}\]
where $\bar{H}_{\nu}^{*2}$ is the tail of $H_{\nu}^{*2},$
so that the inverse Gamma distribution is subexponential in the sense of 
Chistyakov, cf. \cite{chistyakov}. Since it
is moreover infinitely divisible, we deduce by \cite[corollary 1, p.340]{embrechts}
that for all $c>0$\[
\bar{H}_{\nu}^{*c}\left(t\right)\sim\frac{cC_{\nu}}{\nu}t^{-\nu},\,\,\,t\rightarrow +\infty.\]
By Lemma \ref{lem:ultimatelydecreasing} below the density $h_{\nu}^{*c}\left(t\right)$
is ultimately decreasing, so the result follows by application of the
monotone density theorem \cite{Bingham}.
\end{proof}

\begin{lem}
\label{lem:ultimatelydecreasing}
The c-fold convolution of the inverse
Gamma density is ultimately decreasing.
\end{lem}
\begin{proof}
The inverse Gamma density $h_{\nu}\left(t\right)$ is a generalized
Gamma convolution, and so is the convolution power $h_{\nu}^{*c}$, cf. \cite[p. 350]{steutel}.
Since its left extremity is $0,$ we deduce from \cite[prop. 5.5]{steutel}
that it is unimodal, and thus ultimately decreasing.
\end{proof}

\begin{lem}
\label{lem:subordination}
The asymptotic behavior of the c-fold Student
t-convolution is
\[
f_{\nu}^{*c}\left(\mathbf{x}\right)\sim\frac{cA_{\nu}}{\vert \mathbf{x} \vert^{2\nu+d}},\,\, \vert \mathbf{x} \vert \rightarrow\infty.
\]
\end{lem}

\begin{proof}
Since by Lemma \ref{lem:inverseGammaasymptotic} \[
h_{\nu}^{*c}\left(t\right)\sim cC_{\nu}t^{-\nu-1},\]
for any $a$ and $b$ such that $a<cC_{\nu}<b,$ there exists $t_{0}>0$
such that for all $t>t_{0},$\[
at^{-\nu-1}\le h_{\nu}^{*c}\left(t\right)\le bt^{-\nu-1}.
\]
From
\[
\int_{0}^{+\infty}g_{t}\left(\mathbf{x}\right)h_{\nu}^{*c}\left(t\right)dt=\int_{0}^{t_{0}}g_{t}\left(\mathbf{x}\right)h_{\nu}^{*c}\left(t\right)dt+\int_{t_{0}}^{+\infty}g_{t}\left(\mathbf{x}\right)h_{\nu}^{*c}\left(t\right)dt,
\]
it follows that
\begin{eqnarray*}
\int_{0}^{t_{0}}g_{t}\left(\mathbf{x}\right)h_{\nu}^{*c}\left(t\right)dt & + & \int_{t_{0}}^{+\infty}g_{t}\left(\mathbf{x}\right)\frac{a}{t^{\nu+1}}dt  \le \int_{0}^{+\infty}g_{t}\left(\mathbf{x}\right)h_{\nu}^{*c}\left(t\right)dt\\
& \le &\int_{0}^{t_{0}}g_{t}\left(\mathbf{x}\right)h_{\nu}^{*c}\left(t\right)dt+\int_{t_{0}}^{+\infty}g_{t}\left(\mathbf{x}\right)\frac{b}{t^{\nu+1}}dt.
\end{eqnarray*}
But the integral $\int_{0}^{t_{0}}g_{t}\left(\mathbf{x}\right)h_{\nu}^{*c}\left(t\right)dt$
is $o(\vert \mathbf{x} \vert^{-2\nu-d})$ for $\vert \mathbf{x} \vert \rightarrow \infty$ because
\begin{eqnarray*}
\int_{0}^{t_{0}}g_{t}\left(\mathbf{x}\right)h_{\nu}^{*c}\left(t\right)dt &= &\int_{0}^{t_{0}}\frac{1}{\left(4\pi t\right) ^{d/2}}\exp\left(-\frac{\vert \mathbf{x} \vert^{2}}{4t}\right)h_{\nu}^{*c}\left(t\right)dt \\
&\le &\exp\left(-\frac{\vert \mathbf{x} \vert^{2}}{4t_0}\right)\int_{0}^{t_{0}}\frac{1}{\left(4\pi t\right) ^{d/2}}h_{\nu}^{*c}\left(t\right)dt.
\end{eqnarray*}
A simple computation gives
\begin{eqnarray*}
\int_{t_{0}}^{+\infty}g_{t}\left(\mathbf{x}\right)t^{-\nu-1}dt & = & \vert \mathbf{x} \vert^{-2\nu-d}\frac{2^{2\nu}}{\pi^{d/2}}\int_{0}^{\frac{\vert \mathbf{x} \vert^{2}}{4t_{0}}}\exp\left(-u\right)u^{\nu+\frac{d}{2}-1}du\\
& \sim & \vert \mathbf{x} \vert^{-2\nu-d} \frac{2^{2\nu}}{\pi^{d/2}}\Gamma\left(\nu+\frac{d}{2}\right)
\end{eqnarray*}
hence
\begin{eqnarray*}
\limsup_{\vert \mathbf{x} \vert \rightarrow \infty}\vert \mathbf{x}\vert^{2\nu+d}\int_{0}^{+\infty}g_{t}\left(\mathbf{x}\right)h_{\nu}^{*c}\left(t\right)dt
& \le & \limsup_{\vert \mathbf{x} \vert \rightarrow\infty}\vert \mathbf{x}\vert^{2\nu+d}\int_{t_{0}}^{+\infty}g_{t}\left(\mathbf{x}\right)\frac{b}{t^{\nu+1}}dt \\
& = & \frac{b2^{2\nu}}{\pi^{d/2}}\Gamma\left(\nu+\frac{d}{2}\right)
\end{eqnarray*}
and
\begin{eqnarray*}
\liminf_{\vert \mathbf{x} \vert \to \infty}\vert \mathbf{x}\vert^{2\nu+d}\int_{0}^{+\infty}g_{t}\left(\mathbf{x}\right)h_{\nu}^{*c}\left(t\right)dt & \ge & \liminf_{\vert \mathbf{x} \vert \rightarrow \infty}\vert \mathbf{x}\vert^{2\nu+d}\int_{t_{0}}^{+\infty}g_{t}\left(\mathbf{x}\right)\frac{a}{t^{\nu+1}}dt
\\
& = & \frac{a2^{2\nu}}{\pi^{d/2}}\Gamma\left(\nu+\frac{d}{2}\right),
\end{eqnarray*}
so that finally \[
\int_{0}^{+\infty}g_{t}\left(\mathbf{x}\right)h_{\nu}^{*c}\left(t\right)dt\sim c\frac{2^{2\nu}}{\pi^{d/2}}\Gamma\left(\nu+\frac{d}{2}\right)C_{\nu}
\vert \mathbf{x}\vert^{-2\nu-d}
=cA_{d,\nu}\vert \mathbf{x}\vert^{-2\nu-d}.\]
\end{proof}

\begin{lem}
\label{lem:limitcoeff}
For fixed $k$, the coefficients $\alpha_k^{(n)}$ of the Bessel
polynomial of degree $n\ge k$ is increasing in $n$ and 
\[
\lim_{n\rightarrow +\infty} \alpha_{k}^{(n)} = \frac{1}{k!}.
\]
\end{lem}
\begin{proof}
From (\ref{eq:alpha}) we get
\begin{eqnarray*}
\alpha_{k}^{(n)}= \frac{1}{k!}\prod_{j=1}^{k-1}\frac{n-j}{n-\frac{j}{2}} 
\le \frac{1}{k!},
\end{eqnarray*} 
where each of the $(k-1)$ terms of the product
\[
\frac{n-j}{n-\frac{j}{2}} = 1-\frac{\frac{j}{2}}{n-\frac{j}{2}}
\] 
is increasing and converges to $1$ with $n$.
\end{proof}

\begin{lem}
\label{lem:seriesS}
Consider the infinite series $(S)$ equal to
\[
\sum_{k=0}^{+\infty} c_{k} q_{k}(z),
\]
where $q_{k}(z)$ is the $k$'th Bessel polynomial and $c_{k} \in \mathbb{C}$.
Then the three following assertions are equivalent:
\begin{enumerate}
\item
$(S)$ is absolutely convergent for $z=0$,

\item
the sequence $(c_k)$ is absolutely summable,

\item
the series $(S)$ converges absolutely and uniformly on any compact subset of $\mathbb{C}.$
\end{enumerate}
\end{lem}
\begin{proof}
(i) $\Rightarrow$ (ii) since $q_{k}(0)=1.$
\\
(ii) $\Rightarrow$ (iii) since

\begin{eqnarray*}
\vert c_n q_{n}(z)\vert &\le& \vert c_{n} \vert q_{n}(\vert z \vert) \le
\vert c_{n} \vert \sum_{k=0}^{n} \frac{\vert z \vert^{k}}{k !} \le \vert c_{n} \vert
\sum_{k=0}^{+\infty} \frac{\vert z \vert^{k}}{k !} \le \vert c_{n} \vert \exp \left( K \right)
\end{eqnarray*}
for some constant $K.$
The first inequality holds since all the coefficients $\alpha_{k}^{(n)}$ of the Bessel polynomial $q_{n}$ are positive; the second inequality is a consequence of the majorization $\alpha_{k}^{(n)} \le \frac{1}{k!}$ proved in Lemma \ref{lem:limitcoeff}; the third inequality is straightforward, and the last inequality ensues from the assumption that $z$ belongs to a compact subset of $\mathbb{C}.$
Since the sequence $(c_{n})$ is assumed absolutely summable, the absolute and uniform convergence of $(S)$ is a direct consequence  of the above majorization.
\\
(iii) $\Rightarrow$ (i) trivially.
\end{proof}

\section{Conclusion}
In this paper, we have proved two of the conjectures as expressed by Cufaro Petroni in \cite{cufaro}, and disproved a simple version of a third conjecture. We note that these results extend naturally to $d-$dimensional Student t-vectors with correlated components, with density
\[
f_{\nu}(\mathbf{x}) = \frac{A_{d,\nu}}{\vert K \vert}(1+\mathbf{x}^{t}K^{-1}\mathbf{x})^{-(\nu +\frac{d}{2})}
\]
and characteristic function
\[
\varphi(\mathbf{u}) = k_{\nu}(\sqrt{\mathbf{u}^{t}K\mathbf{u}})
\]
where $K$ is a symmetric and positive definite matrix; Theorem \ref{thm:thm1} holds unchanged and the asymptotic result of Theorem \ref{thm:cstepStudent} still holds by replacing constant $A_{d,\nu}$ by $\frac{A_{d,\nu}}{\vert K\vert^{\frac{1}{2}}}$.

\section*{References}
\bibliography{cufarobib}
\bibliographystyle{jphysicsB}
\end{document}